\documentclass[a4paper,
12pt]{article}
\usepackage{epsfig}
\usepackage{url}
\usepackage{graphicx}
\usepackage{geometry}
\usepackage{hyperref}
\usepackage{xcolor}
\usepackage{url} 
\usepackage{amsmath}
\usepackage{amsthm}
\usepackage{amssymb} 
\usepackage{latexsym}
\usepackage{authblk}
\DeclareMathOperator*{\argmin}{\arg\!\min}
\begin{document}
\newtheorem{theorem}{Theorem}[section]
\newtheorem{corollary}[theorem]{Corollary}
\newtheorem{lemma}[theorem]{Lemma}
\newtheorem{proposition}[theorem]{Proposition}
\newtheorem{remark}[theorem]{Remark}
\newtheorem{conjecture}[theorem]{Conjecture}
\newtheorem{assumption}[theorem]{Assumption}
\newcommand{\ind}{1\hspace{-2.5mm}{1}}
\newcommand{\amax}{a_{\text{max}}}
\newcommand{\Bin}{\text{Bin}}

\newcommand{\red}[1]{\textcolor{red}{#1}}
\newcommand{\til}[1]{\widetilde{#1}}

\title{Who is the infector? General multi-type epidemics and real-time susceptibility processes}
\author{Tom Britton$^1$, Ka Yin Leung$^1$ and Pieter Trapman$^1$}
\maketitle

\begin{abstract}
We couple a  multi-type stochastic epidemic process with a directed random graph, where edges have random lengths.
This random graph representation is used to characterise the fractions of individuals infected by the different types of vertices among all infected individuals in the large population limit.   
For this characterisation we rely on theory of multi-type real-time branching processes.
We identify a special case of the two-type model, in which the fraction of individuals of a certain type infected by individuals of the same type, is maximised among all two-type epidemics approximated by branching processes with the same mean offspring matrix. 
\end{abstract}


\noindent2010 Math.\ Subj. Classification:
Primary 92D30; 60K35; 
Secondary 05C80; 60J80

\footnotetext[1]{Department of Mathematics, Stockholm University, 106 91 Stockholm, Sweden.\\ Email: \{tom.britton, kayin.leung, ptrapman\}@math.su.se}
\section{Introduction}
Mathematical models have proven to be successful in understanding infectious disease dynamics. Often, the focus is on (controlling) the beginning of an epidemic. In many of those models branching process approximations and the concept of the basic reproduction number $R_0$ (which corresponds to the offspring mean of the approximating branching process) play an important role. In the current paper we consider  the entire epidemic outbreak instead. We do so in a setting where there are multiple types of infected individuals. Suppose that a large epidemic outbreak has taken place in the population. Then a certain fraction of the population will have become infected. The question that we concern ourselves with is: what fraction of the infected population of a certain type $j$ was infected by individuals of type $i$? Here $i,j=1,\ldots K$, where $K$ is the number of types in the population. In other words, who is the infector?

This question is not so straightforward to answer. Timing of events plays an essential role. We approach this question using an epidemic graph construction, which is used as a tool in proving the two main theorems of this paper. But, the construction proves to be interesting in itself, and a substantial part of this paper is devoted to this construction. Using the graph representation of the epidemic, we consider susceptibility processes. Susceptibility \emph{sets} were introduced in infectious disease modelling by Ball and co-authors \cite{Ball01,Ball02}. The susceptibility set of a vertex $v$ consists of all vertices $u$ in the vertex set from which there is a path from $u$ to $v$ in the (restricted) epidemic random graph. The epidemic process and epidemic random graph can be coupled in such a way that $u$ is in the susceptibility set of $v$ if and only if $v$ is infected during the epidemic conditioned on $u$ being initially infectious. Susceptibility sets have proven to be important tools in proving results concerning e.g.\ the final size of different epidemic models \cite{Ball02,Ball09} (see also \cite{Ande99}). However, for the research question in this paper we need to consider the susceptibility process instead, i.e.\ we need to take timing into account. 

In this paper, we also consider a multi-type (backward) branching process. This branching process is constructed in such a way that the distribution of the tree-like graph corresponding to it, is the same as the susceptibility process (up to a certain time). The coupling between the susceptibility process and the branching process enables us to prove the main results of this paper formulated in Theorems~\ref{firstmain} and \ref{secondmain}. Using existing theory for branching processes we find an answer in Theorem~\ref{firstmain} for the question ``who is the infector?'' by means of an expression for the expected fraction $\rho_{ij}$ of infected individuals of type $j$ that are infected by individuals of type $i$, conditioned on that there is a large outbreak in a population (as the population size tends to infinity), $i,j=1,\ldots K$. 
In general, this expression remains rather implicit. However, for a special class of models, we are able to obtain upper and lower bounds for the quantities $\rho_{ij}$ of interest (Theorem~\ref{secondmain}) if we keep the fractions of individuals of the different types and the expected number of infectious contacts between different types of individuals fixed. As the $\rho_{ij}$ of Theorem~\ref{secondmain} are rather implicit, these bounds allow us to gain more insights in the importance of different types of infected individuals in the transmission dynamics in the population. 

The class of models that allows for identifying upper and lower bounds are the topic of interest of our twin paper \cite{Leun18}. This paper \cite{Leun18} is motivated by infectious diseases, such as influenza and chlamydia, for which we can categorise infected individuals as symptomatic or asymptomatic (showing no apparent signs of the disease) giving rise to two types of infected individuals. Asymptomatically infected individuals are generally hard to detect by public health authorities. Therefore, we would like to gain insights in their role in the transmission process and determine whether asymptomatically infected individuals often play the role of the infector. 

The structure of this paper is as follows. In Section~\ref{secmod} we introduce the model, the notation and the two main theorems of the paper. Next, in Section~\ref{sec:construction}, we discuss the construction that enables us to prove the desired results. The proofs are then presented in Section~\ref{sec:infector}. We end with a short discussion in Section~\ref{sec:discussion}.

\section{Model, notation and main results}
\label{secmod}

\subsection{Model and notation}\label{sec:model_notation}
We denote the number of elements in a set $\mathcal{A}$ by $|\mathcal{A}|$. For $k \in \mathbb{N}$ we use the notation $[k] = \{1,2, \cdots, k\}$. We use $\mathbb{N}$ for the strictly positive integers and $\mathbb{N}_0 = \mathbb{N} \cup \{0\}$ for the non-negative integers. 

Unless specified otherwise, limits are for population size $n \to \infty$. We say that an event happens with high probability (w.h.p.) if the probability of the event converges to 1 as $n \to \infty$. We adhere to the usual order notation, i.e.\ $f =O(g)$ means that $\limsup_{n \to \infty} |f(n)/g(n)| < \infty$ and $f =o(g)$ means that $\lim_{n \to \infty} |f(n)/g(n)| =0$. In addition, for a sequence of random variables $\{X^{(n)}; n\in \mathbb{N}\}$, we write $X^{(n)} = O_p(g)$ if $|X^{(n)}/g(n)|$ is bounded in probability and $X^{(n)}=o_p(g)$ if $X^{(n)}/g(n) \to 0$ in probability. See \cite[Section 1.2]{Jans11} for a discussion of this notation.  

We consider a population of $n$ individuals where $V^{(n)}$ denotes the set of all individuals. For some of our results we consider a sequence of models in growing populations, i.e.\ for $n \to \infty$. If no confusion is possible we write $V= V^{(n)}$. We assume that there are $K$ types of individuals. For $i \in [K]$, let $V_i = V_i^{(n)}$ be the set of vertices of type $i$ and $n_i = |V_i|$ the number of vertices of type $i$. 


Within this sequence of populations we consider the spread of an infectious disease in which individuals are either susceptible or infected. When individual $v \in V_i^{(n)}$ is infected it makes contacts to different individuals of type $j$ ($i,j \in [K]$) according to a point process $\xi_v^j = \{\xi_v^j(t);t\geq 0\}$ on the interval $[0,\infty)$. If there are more than $n_j$ points in the process, then only the first $n_j$ points represent contacts. The time parameter in the definition of $\xi_v^j$ represents the time since infection of individual $v$. The processes $\xi_v^1$, $\xi_v^2$, $\cdots$ and  $\xi_v^K$ may be dependent, and their joint distribution may depend on the type $i$ of $v$. However, the point processes associated to different individuals are independent, i.e.\ the vectors of processes $\{(\xi_v^1, \cdots \xi_v^K);v \in V\}$ are independent. For convenience we introduce the vector of stochastic process $(\xi_{i1}, \xi_{i2},\cdots, \xi_{iK})$, which is distributed as $(\xi_v^1, \cdots \xi_v^K)$ for $v \in V_i$, $i \in [K]$. We make the following assumptions.

\begin{assumption}\label{nipiass}
$n^{-1} n_i \to p_i>0$  for all $i \in [K]$. 
Furthermore, 
\begin{equation*}
\max_{i,j\in [K]} \left|\frac{p_i}{p_j}-\frac{n_i}{n_j}\right|=O(1/n).
\end{equation*}
\end{assumption}

\begin{assumption}\label{finass}
For all $i \in [K]$, the distribution of $(\xi_{i1}, \xi_{i2},\cdots, \xi_{iK})$ is independent of population size $n$ and  for all $i,j \in [K]$,
\begin{equation}\label{eq:mij}
m_{ij}=\mathbb{E}[\xi_{ij}(\infty)]<\infty,
\end{equation}
i.e.\ the expected number $m_{ij}$ of contacts that an individual of type $i$ makes with individuals of type $j$  is finite, for all $i,j \in [K]$.  Furthermore, we assume that $|(\xi_{i1}, \xi_{i2},\cdots, \xi_{iK})|$ does a.s.\ only contain jumps of size 1 and its distribution has no atoms. 
\end{assumption}

\begin{assumption}\label{boundass}
There exists a constant $\kappa <1$ such that for all $i \in [K]$, $$\mathbb{P}\left(\max_{v \in V} \xi_v^i(\infty) < n^{\kappa}\right) \to 1 \qquad \mbox{as} \qquad n \to \infty.$$
\end{assumption}

For future reference, we let $M=\{m_{ij}\}_{i,j\in[K]}$ denote the matrix with the $m_{ij}$ (as defined in Assumption \ref{finass}) as elements.
Assumption \ref{finass} guarantees that w.h.p.\ all non-trivial paths in the epidemic graph defined below have different lengths. 
Note that Assumption \ref{boundass} is easily met, e.g.\ if for all $v \in V$ and $i \in [K]$ there exists $\epsilon>0$ and $\ell_0 \in (0,\infty)$, for which $\mathbb{P}[\xi_v^i(\infty)>\ell] < \ell^{-(1+\epsilon)}$ for all $\ell>\ell_0$, then for $\kappa \in (1/(1+\epsilon),1)$
\begin{align*}
\mathbb{P}\left(\max_{v \in V} \xi_v^i(\infty) < n^{\kappa}\right) &=  \prod_{v \in V} [1-\mathbb{P}\left( \xi_v^i(\infty) > n^{\kappa}\right)] \\ 
&\geq 1 - \sum_{v \in V}\mathbb{P}\left( \xi_v^i(\infty) > n^{\kappa}\right) \\ 
&\geq 1 - n \max_{j \in [K]}\mathbb{P}\left( \xi_{ji}(\infty) > n^{\kappa}\right)\\
& >1 -n \times n^{-(1+\epsilon)\kappa} \\
& \to 1.
\end{align*}

At the points of $\xi_v^j$ ($j \in [K]$), $v$ contacts an individual from $V_j$. The individuals that are contacted are uniformly chosen without replacement. If $v$ is of type $j$, we allow for $v$ to be among the contacted individuals in $\xi_v^j$. If the individual that is contacted is still susceptible at the time of the contact then it becomes infected. 

Note that we may assign the point processes $\{\xi_v^j;v \in V, j \in [K]\}$  and decorate the points with the labels of the individuals these points represent contacts to, already before the epidemic starts. This allows us to create a new set of random variables $\{\eta(u,v);u,v \in V\}$. Assume that individual $v$ is of type $j \in [K]$. If there is a point in $\xi_u^j$ with label $v$, then $\eta(u,v)$ takes the value of this point. If there is no such point in $\xi_u^j$, set $\eta(u,v)= \infty$. Observe that the distribution of $\{\eta(u,v);u,v \in V\}$ depends on $n$ in this construction, because the probability that $v$ is chosen as a label is decreasing in $n$. Note that there is a broad class of models that satisfy  Assumptions~\ref{nipiass}-\ref{boundass}, see Remark~\ref{rk:SEIR} for an example of a specific class of epidemic models.

\begin{figure}
\centering
\includegraphics[scale=0.7]{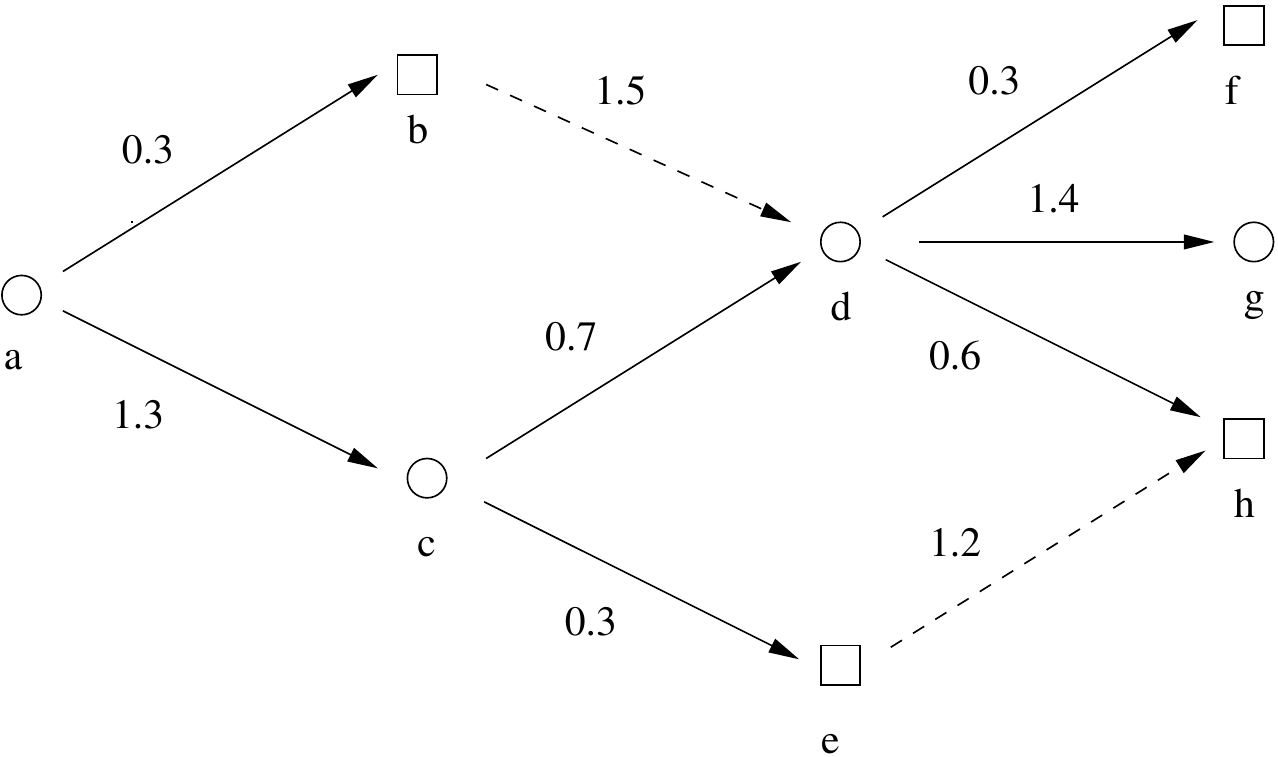}
\caption{An example of $G'=(V,E')$ with $n=8$ and $K=2$. The vertices are labeled $a,b,\cdots,h$. Vertices of type 1 are represented by circles  and vertices of type 2 by boxes. Contacts made by  vertices of type 1 are represented by solid directed edges and contacts made by  vertices of type  2 by dashed directed edges. The numeric values next to the edges reflect the time since infection $\eta(u,v)$ of the tail $u$ of the edge until a contact with the head $v$ takes place, $u,v \in \{a,b,\cdots,h\}$.}
\label{fig:forward}
\end{figure}

This construction provides us with a graph representation of the population and the epidemic on it. We construct the weighted random graph $G=(V,E)$ as follows. The edge set $E$ consists of all directed pairs $(u,v) \in V \times V$, with $u \neq v$. For edge $(u,v) \in V \times V$, we say that $u$ is the tail of $(u,v)$ and $v$ is its head. The weight of edge $(u,v)$ is given by $\eta(u,v)$ for all $(u,v) \in E$. For some of our arguments we restrict to the weighted edge set $E' \subset E$ of all edges $(u,v) \in E$ with finite weight, i.e.\ $(u,v) \in E'$ if and only if  $\eta(u,v)< \infty$. The corresponding random graph is denoted by $G'=(V, E')$. For $i,j \in [K]$,  $u \in V_i$ and $v \in V_j$, we also introduce the random variable $\eta_{ij}$, which is distributed as $\eta(u,v)$ conditioned on $\eta(u,v)<\infty$ (i.e.\ conditioned on $(u,v) \in E'$).
Observe that 
\begin{equation}\label{disteta}
\begin{aligned}
\mathbb{P}(\eta_{ij} \leq t) &= \mathbb{P}(\eta(u, v) \leq t|\eta(u, v) < \infty) \\
&= \frac{\mathbb{P}(\eta(u,v) \leq t)}{\mathbb{P}(\eta(u,v) < \infty)} = \frac{\frac{1}{n_j}\mathbb{E}[\xi_{ij}(t)]}{\frac{1}{n_j}\mathbb{E}[\xi_{ij}(\infty)]}= \frac{\mathbb{E}[\xi_{ij}(t)]}{m_{ij}},
\end{aligned}
\end{equation}
which is independent of $n$. We mainly consider the random graph $G'=(V,E')$. See Figure \ref{fig:forward} for an example of $G'$.

An ordered set of distinct vertices $\pi =(v_1, v_2, \cdots, v_m)$ is a path in $G'$ if $(v_i,v_{i+1}) \in E'$ for all integers $i \in [m-1]$. With some abuse of terminology, we sometimes refer to the set of edges connecting the vertices of $\pi$ as the path $\pi$ and speak of a path in $E'$. The length of a path $\pi =(v_1, v_2, \cdots, v_m)$ is $\ell(\pi) = \sum_{i=1}^{m-1} \eta(v_i,v_{i+1})$. Assumption \ref{finass} implies that all non-zero paths in $E'$ have different lengths with probability 1. Let $\Pi_{uv}$ be the set of all paths from $u$ to $v$ in $E'$ and define the (quasi) distance from $u$ to $v$ as $d(u,v)=\min_{\pi \in \Pi_{uv}}\ell(\pi)$. As an example, in Figure \ref{fig:forward} the distance from $a$ to $d$ is given by $$d(a,d) =\min\left(\eta(a,b)+\eta(b,d), \eta(a,c)+\eta(c,d)\right)=1.8.$$
In general, $d(u,v)\neq d(v,u)$ if $u\neq v$ since the graph $G'$ is directed.   Therefore $d$ is actually a quasi-distance. We say that $d(u,v)=\infty$ if $\Pi_{uv} = \emptyset$. Furthermore, $d(v,v)=0$ for all $v\in V$. 

Next we formulate an ``irreducibility'' assumption:
\begin{assumption}
\label{assirred}
For every $i,j\in[K]$, there is w.h.p.\ a path from a vertex in $V_i$ to a vertex in $V_j$ in $E'$.
\end{assumption}



An epidemic process is reproduced from $G'$ as follows. 
Let $V_{\text{init}}=V_{\text{init}}^{(n)}$ be the set of vertices that are initially infected. This set may be predetermined or randomly selected and satisfies the following assumption.
\begin{assumption}
\label{assvinit}
The set $V_{\text{init}}^{(n)}$ of vertices that are initially infected satisfies $|V_{\text{init}}^{(n)}|=O_p(1)$.
\end{assumption} 

We set $\sigma_v$ to be the time between the start of an epidemic until $v$ becomes infected, $v \in V$, i.e.\ $
\sigma_v = \inf \{d(u,v);u \in V_{\text{init}}\}$. In particular, $\sigma_v =0$ for $v \in V_{\text{init}}$.
So, suppose we let $V_{\text{init}}=\{a\}$ in Figure \ref{fig:forward}, then $\sigma_a=0$, $\sigma_b=0.3$, $\sigma_c = 1.3$,  $\sigma_d = 1.8$ etc.
Moreover, if $\sigma_v= \infty$, then $v$ will never become infected. Note that $G'$ contains some redundant information regarding the epidemic: (i) there are edges with finite weights with infection time of the tail being $\infty$, i.e.\ the tails of those edges will never get infected, (ii) if $\sigma_v < \sigma_u + \eta(u,v)$, then the edge $(u,v)$ represents a contact between two already infected individuals. 

\begin{remark}[The SEIR epidemic model]\label{rk:SEIR}
The framework of this section includes the SEIR (Susceptible $\to$ Exposed $\to$ Infectious $\to$ Recovered) epidemics as follows. In the SEIR framework, exposed, infectious and recovered individuals are all counted as infected. Assign to every vertex $v \in V$ a random latent period $L_v$ and a random infectious period $\iota_v$, which might be dependent on $L_v$. The vectors $\{(L_v,\iota_v);v \in V\}$ are independent and their distribution functions only depend on the type of the vertex. Suppose that $v$ is of type $i$. Conditioned on $(L_v,\iota_v)$, let the processes $\{\hat{\xi}_v^j; j \in [K]\}$ be independent homogeneous Poisson processes on the interval $(L_v, L_v + \iota_v)$ with intensity $p_j \lambda_{ij}$. At the points of this point process, $v$ makes contacts to  vertices in $V_j$ chosen uniformly with replacement. By keeping only the points in $\{\hat{\xi}_v^j; j \in [K]\}$ whose label did not appear before in this process, we obtain $\{\xi_v^j; j \in [K]\}$. If $L_v$ and $\iota_v$  are exponentially distributed, then the SEIR epidemic is a Markov process that is often referred to as the Markov SEIR epidemic.

At time $t$, a vertex $v$ is susceptible if $t< \sigma_v$, Exposed if $t \in [\sigma_v, \sigma_v + L_v)$, while it is Infectious if $t \in [\sigma_v + L_v, \sigma_v + L_v+ \iota_v)$. Finally $v$ is Recovered if $t \geq \sigma_v + L_v+ \iota_v$.

If $\mathbb{P}(L_v =0)=1$, then we are in the so-called SIR (Susceptible $\to$ Infectious $\to$ Recovered) epidemics framework. In addition, if $\iota_v$ is exponentially distributed then the process is called a Markov SIR epidemic. 
\end{remark}

Finally, to conclude this section, we introduce the basic reproduction number $R_0$. As this is possibly the most studied quantity in mathematical modelling of the spread of infectious diseases, no epidemic modelling paper would be complete without at least mentioning $R_0$. In a single-type epidemic model $R_0$ is defined as the expected number of infectious contacts made by a newly infected individual in an otherwise susceptible population. The multi-type equivalent of $R_0$ is given by the dominant eigenvalue of the $K \times K$ matrix $M=\{m_{ij}\}_{i,j\in[K]}$ (see Assumption~\ref{finass}). $R_0$ is always real and strictly positive \cite[Chapter 4]{Jage75},\cite[Chapter 7]{Diek12}. We say that the epidemic process is supercritical (resp.\ critical, resp.\ subcritical) if $R_0>1$ (resp.\ $R_0=1$, resp.\ $R_0<1$). As $n \to \infty$, the epidemic becomes large  with positive probability if and only if the process is supercritical \cite{Diek12}. In the remainder of the paper we assume that $R_0>1$.


\subsection{Main results}
Because we assume that $R_0>1$, a large outbreak occurs with positive probability and at the end of such a large outbreak, a fraction of the infected individuals is of type $j$ w.h.p. However, the question that we are interested in is: what fraction of those individuals was infected by individuals of type $i$, for $i,j\in[K]$?
So, if we denote set of infected individuals of type $j$ by $\mathcal{I}_j$, and the subset of those individuals, which are infected by individuals of type $i$  by $\mathcal{I}_j^i$, then we are interested in $|\mathcal{I}_j^i|/|\mathcal{I}_j|$. 
For example, in Figure \ref{fig:forward}, if $a$ represents the initially infectious individual,
then there are two individuals of type 1 that are infected by an individual of type i (namely the individuals represented by vertices $c$ and $g$) and one individual of type 1 that is infected by an individual of type 2 (namely the individual represented by vertex $d$).
So, the fraction of individuals infected by an individual of type 1 among all infected  individuals of type 1 is $2/3$ (we exclude vertex $a$ because it was initially infected, and not infected by another individual in the population).
The  question posed above leads to the main results that are formulated in Theorems~\ref{firstmain} and~\ref{secondmain} below. For our results, we need to define a multi-type branching process $\{\mathcal{Z}(t); t \geq 0\}$. The branching process is defined as follows. Let $Z^i(t)$ be the number of particles in the branching process at time $t$ if the process is started by a single  particle of type $i$, $i \in [K]$. For $i,j \in [K]$, particles of type $j$ in $\mathcal{Z}(t)$ give birth to particles of type $i$ according to a Poisson process with intensity $\frac{p_i}{p_j} \mathbb{E}[\xi_{ij}(da)]$, where $a$ is the age of the particle of type $j$. All of these Poisson processes are independent.

The definition of the branching process $\mathcal{Z}(t)$ is such that we can apply theory from $\cite{Iksa15}$. In particular, there exists a Malthusian parameter $\alpha>0$ and a random variable $W^i$, such that 
\begin{equation*}
e^{-\alpha t}Z^i(t) \to W^i \qquad \mbox{ a.s.\ as $t \to \infty$}
\end{equation*} 
and 
\begin{equation*}
\mathbb{P}(W^i \in (0, \infty))=1 -\mathbb{P}(W^i =0) = \mathbb{P}(Z^i(t) \to \infty).
\end{equation*} 
Let $W^i(r)$ for $r \in \mathbb{N}$ be independent copies of $W^i$. We now state our first main result. 

\begin{theorem}
\label{firstmain}
Conditioned on the occurrence of a large outbreak,
the fraction of infected individuals of type $j$ that are infected by an individual of type $i$ during an outbreak in a population of size $n$ converges in probability to $\rho_{ij}$. Here $\sum_{i=1}^K\rho_{ij}=1$, and
\begin{equation}
\label{mainequ}
\rho_{ij}=\frac{1}{\mathbb{P}(Z^j(t) \to \infty)}
\mathbb{E}\left[\frac{\sum_{r=1}^{X_{ij}}e^{-\alpha \tau_{ijr}}  {W}^i(r)}{\sum_{k =1}^K \sum_{r=1}^{X_{kj}} e^{-\alpha \tau_{kjr}}  {W}^k(r)} \ind\left(\sum_{k =1}^K \sum_{r=1}^{X_{kj}}  {W}^k(r) >0\right)
\right],
\end{equation}
where, for $k \in [K]$ and $r \in \mathbb{N}$, the random variables $\tau_{kjr}$ are independent with distribution function $\mathbb{P}(\tau_{kjr} <a) = {\mathbb{E}[\xi_{kj}(a)]}/{m_{kj}}$ and $X_{kj}$ is Poisson distributed with expectation $\frac{p_k}{p_j}m_{kj}$.
\end{theorem}

Note that, in general, it is hard to give a more explicit expression for $\rho_{ij}$ than~\eqref{mainequ}. Often, there is no explicit description of the distribution of ${W}^k(r)$. We are able to obtain bounds for~\eqref{mainequ} for the important special case discussed in \cite{Leun18}, leading to Theorem~\ref{secondmain}. 

The model of \cite{Leun18} is as follows. We consider $K=2$, and $(\xi_v^1,\xi_v^2)$, $v \in V$, obtained from a single marked point process $\xi_v$. In this process $\xi_v$ the points get independently mark 1 with probability $p_1$ and mark 2 otherwise. Then the process $\xi_v^1$ consists of the points with mark 1, while $\xi_v^2$ consists of the points with mark 2. By construction of the process, the probability that an ultimately infected vertex of type $i$ is infected by a  vertex of type 1 is the same for $i=1$ and $i=2$, i.e.\ $\rho_{11}=\rho_{12}=\rho_1$. With some abuse of notation we write $\xi_i$ for a point process distributed as $\xi_v$ for $v \in V_i$, $i \in [2]$. Furthermore, note that $m_{ij}=\mathbb E[\xi_{ij}(\infty)]=p_j\mathbb E[\xi_i(\infty)]$, i.e.\ we can write $m_{ij}=p_j \tilde m_i$ with $\tilde m_i=\mathbb E[\xi_i(\infty)]$. Here $\tilde m_i$ can be interpreted as the expected number of secondary cases caused by a newly infected individual of type $i$ in an otherwise susceptible population (one can think of the $\tilde m_i$ as the type-specific reproduction numbers). Then the basic reproduction number is $R_0=p_1\tilde m_1+p_2\tilde m_2$, with $p_2=1-p_1$. Indeed, a newly infected individual is of type $i$ with probability $p_i$ and the expected number of secondary cases it produces is $\tilde m_i$, $i=1,2$.

For this model we can compute $\rho_1^-$ and $\rho_1^+$, the minimum and maximum fraction of the infected vertices that are infected by type 1 vertices, if the matrix $M$ and $p_i$ are held fixed for $i \in [K]$. We let $\tilde{q}_1$ be the smallest positive solution in $(0,1]$ of 
\begin{equation}\label{qsequa}
x= e^{-p_1 \tilde m_{1}(1-x)}
\end{equation} 
and $\tilde q_2$ the smallest positive solution in $(0,1]$ of
\begin{equation}\label{qaequa}
x= e^{-p_2 \tilde m_{2}(1-x)}.
\end{equation} 
Furthermore, we let $q$ be the unique solution in $(0,1)$ of 
\begin{equation}\label{qequa}
x  =  e^{-(1-x)\left(p_1 \tilde m_{1}+ p_2 \tilde m_{2}\right)}=e^{-(1-x)R_0}.
\end{equation}
Note that we assume that $R_0=p_1\tilde m_{1}+ p_2\tilde m_{2}>1$, so the unique solution $q\in(0,1)$ exists. We can interpret $1-\tilde q_1$ (resp.\ $1-\tilde q_2$) as the final fraction of the population that ultimately gets infected when only individuals of type 1 (resp.\ type 2) are able to transmit, conditional on a large outbreak. Furthermore, $1-q$ can be interpreted as the final fraction of the population that ultimately gets infected, conditional on a large outbreak (or conversely, $q$ is the fraction of the population that remains susceptible throughout the epidemic). 

\begin{theorem}
\label{secondmain}
Consider the two-type model described above. In the limit as population size $n\to\infty$ and for fixed $p_1$, $\tilde{m}_1$ and $\tilde{m}_2$, the fraction of ultimately infected vertices that is infected by type 1 vertices is bounded from above by 
\begin{equation*}
\rho_1^+ = \left( 1 -  \frac{p_1 \tilde m_1 (\tilde{q}_1 + q)}{2} \right) \frac{\tilde{q}_1-q}{1-q}.
\end{equation*}
\end{theorem}

By interchanging the role of the types $1$ and $2$, we also obtain the lower bound $\rho_1^-$. Indeed, note that $\rho_1^-=1-\rho_2^+$, where
\begin{equation*}
\rho_2^+ = \left( 1 - \frac{p_2 \tilde m_2 (\tilde{q}_2 + q)}{2} \right) \frac{\tilde{q}_2-q}{1-q}.
\end{equation*}
In other words, for any  point process $\{(\xi_v^1,\xi_v^2); v\in V\}$ that satisfies the assumptions of Section~\ref{sec:model_notation} and that can be obtained from independently marking points of a one-dimensional point process, the fraction $\rho_1$ of infected individuals that are infected by individuals of type 1 is bounded by $\rho_1^-$ and $\rho_1^+$, i.e.\ $\rho_1^-\leq \rho_1 \leq \rho_1^+$. 

\begin{remark}
In Section~\ref{sec:proofsecondmain} we consider a more general setting than the one in Theorem~\ref{secondmain}. Instead of assuming a single marked process $\xi_i$, one can consider general distributions $(\xi_{i1}, \xi_{i2})$ and obtain bounds~\eqref{finaleq} for $\rho_{21}^-$ and~\eqref{eq:rho11_final} for $\rho_{11}^+$ (and, by interchanging the roles of types 1 and 2, bounds $\rho_{12}^-$ and $\rho_{22}^+$). As this is somewhat more involved, we choose to present the bounds in the form of Theorem~\ref{secondmain} instead.
\end{remark}

\section{Susceptibility process, backward branching process and the coupling}\label{sec:construction}

Throughout we assume that all random variables and processes are defined on a suitable rich enough probability space, which we do not specify.

\subsection{The susceptibility process}\label{sec:susceptibility}
In this subsection we use the idea of susceptibility sets \cite{Ball01,Ball02,Ball09}, and construct this set through a stochastic process: the susceptibility process.

We define the susceptibility process $\{\mathcal{S}^{(n)}_v(t);t \geq 0\}$ as 
\begin{equation*}
\mathcal{S}^{(n)}_v(t) = \{u \in V; d(u,v) \leq  t\}.
\end{equation*}  
Note that $\mathcal{S}^{(n)}_v(t)$ is non-decreasing in $t$. The susceptibility set of vertex $v$ is defined as $\mathcal{S}^{(n)}_v = \lim_{t \to \infty} \mathcal{S}^{(n)}_v(t)$, i.e.\ the susceptibility set $S_v^{(n)}$ of $v$ consists of all vertices $u\in V$ for which there is a path from $u$ to $v$ in $G'$. As an illustration, in Figure \ref{fig:forward}, the susceptibility set of vertex $f$  is given by  $\mathcal{S}^{(n)}_f = \{a,b,c,d\}$, while $\mathcal{S}^{(n)}_f(t=1.1) = \{c,d\}$. Note that $\mathcal{S}^{(n)}_v \cap V_ {\text{init}} = \emptyset$ if and only if $v$ remains uninfected throughout the epidemic, i.e.\ if and only if there is no path in $G'$ from $V_ {\text{init}}$ to $v$. Also note that the susceptibility set may contain vertices of different types. 

Define for $a \geq 0$ and $j \in [K]$, 
\begin{equation*}
\mathcal{S}^{(n)}_v(t;a,j) = \{u \in \mathcal{S}^{(n)}_v(t)\cap V_j; d(u,v) > t- a\}.
\end{equation*}
That is, $\mathcal{S}_v(t;a,j)$  consists of the  vertices of  type $j$ in $\mathcal{S}^{(n)}_v(t)$, that are not yet part of $\mathcal{S}^{(n)}_v(t-a)$.  

Let $v \in V \setminus V_{\text{init}}$ be a randomly chosen vertex. We derive the susceptibility process $\{\mathcal{S}^{(n)}_v(t);t \in (0,t_*)\}$ by constructing part of the random graph $G'$ around vertex $v$ by means of an exploration process $\{\hat{G}(\ell)\}= \{\hat{G}(\ell); \ell \in \mathbb{N}_0\}$ in which vertices in the susceptibility process are explored one at a time. We note that $\{\hat{G}(\ell)\}$ depends on $v$. Here $t^*= t^*(n)$ is a given time (we defer the specification of $t^*$ until~\eqref{eq:t*} below). The process $\{\hat{G}(\ell)\}$ allows us to couple $\{\mathcal{S}^{(n)}_v(t)\}$ with an appropriate branching process. In this way, we can make the coupling between the susceptibility process and the backward branching process that is needed to prove Theorems~\ref{firstmain} and~\ref{secondmain} in Section \ref{sec:infector}.

Before we define the exploration process $\{\hat{G}(\ell)\}$ around vertex $v$, we introduce some additional variables and terminology. $\hat{G}(\ell)$ is a 4-tuple: 
\begin{equation*}
\{\hat{G}(\ell)= (\hat{V}^{a}(\ell),\hat{V}^{p}(\ell),\hat{V}^{e}(\ell),\hat{E}(\ell));\ell \in \mathbb{N}_0\}.
\end{equation*}
Here $\hat E(\ell)$ denotes the edge set of $\hat{G}(\ell)$. Vertices in $\hat{G}(\ell)$ can be `active', `passive', or `explored'. The sets of these vertices are denoted by $\hat{V}^{a}(\ell)$, $\hat{V}^{p}(\ell)$, and $\hat{V}^{e}(\ell)$, respectively. The sets $\hat{V}^{a}(\ell)$, $\hat{V}^{p}(\ell)$, $\hat{V}^{e}(\ell)$, and $\hat{E}(\ell)$ are defined in the construction below (where also their names will become apparent). Finally, before we explain the construction, we mention that throughout the process we may `flag' the process (see step 4). This flagging plays a role when coupling the exploration process with a branching process to represent the susceptibility process. The construction is as follows. 
\begin{itemize}
\item[1)] 
To set the initial conditions of the construction, let $\hat{V}^a(0) = v$, $\hat{V}^e(0) = \emptyset$, and $\hat{E}(0)$ be the set of all edges in $E'$ for which $v$ is the tail, whereas $\hat{V}^p(0)$ is the set of all heads of edges in $\hat{E}(0)$ that are in $V\setminus v$.

\item[2)] 
For $\ell \in \mathbb{N}_0$, assume that $\hat{V}^a(\ell)\neq \emptyset$ and that there exists a vertex in $\hat{V}^a(\ell)$ from which there is a path in $\hat{G}(\ell)$ to $v$ of length at most $t^*$. In step $\ell+1$ pick (according to some rule) one of the vertices from $\hat{V}^a(\ell)$, from which there is a path in $\hat{G}(\ell)$ to $v$ of length at most $t^*$. Say that this vertex is $v' \in V_j$. Move $v'$ to the set of explored vertices, i.e.\ $\hat{V}^e(\ell+1) = \hat{V}^e(\ell) \cup v'$. Assign to $v'$ independently a binomial random number $x(v';i)$ with parameters $n_i$ and $m_{ij}/n_j$.

\item[3)] 

The remainder of step $\ell+1$ is split up  in $\sum_{i=1}^K x(v';i)$ sub-steps as follows. We introduce
\begin{equation*}
\left\{{G}^*(\ell,\ell');\ell \in \mathbb{N}_0, \ell' \in \{0,1,\cdots, \sum_{i=1}^K x(v';i)\}\right\},
\end{equation*}
where
\begin{equation*}
{G}^*(\ell,\ell')= ({V}^{*,a}(\ell,\ell'),\hat{V}^{*,p}(\ell,\ell'),\hat{V}^{*,e}(\ell,\ell'),\hat{E}^*(\ell,\ell')).
\end{equation*}
Furthermore, set ${G}^*(\ell,0)=\hat{G}(\ell)$. Next, let $\ell' \in \left(\sum_{i'=1}^{i-1} x(v';i'),\sum_{i'=1}^i x(v';i')\right]$, where the empty sum is $0$.
\item[4)]
In the $\ell'$-th sub-step of step $\ell+1$ pick uniformly a vertex from $V_i$.
If we pick a vertex we have picked in one of the $\ell'-1$ sub-steps before we say that the exploration process is \textit{flagged}. 
In that case we choose new vertices from $V_i$ until we obtain a vertex that was not chosen in the previous $\ell'-1$ sub-steps. This step is equivalent to picking the vertices without replacement.
Say that the vertex that is picked is $v''$. 
\begin{itemize}
\item [4a)] If $v'' \in \hat{V}^{a}(\ell) \cup \hat{V}^{e}(\ell)$, then nothing changes in the exploration graph, i.e.\ ${G}^*(\ell+1,\ell')= {G}^*(\ell+1,\ell'-1)$. This is because if $v'' \in \hat{V}^{a}(\ell) \cup \hat{V}^{e}(\ell)$, then we already have explored the edges with tail $v''$.  

\item [4b)] If $v'' \in V \setminus \left(\hat{V}^{a}(\ell) \cup \hat{V}^{e}(\ell)\right)$, assign to $v''$ the vector of point processes $(\xi_{v''}^{j'}, j' \in [K])$. The distribution of $(\xi_{v''}^{j'}, j' \in [K])$ is equal to the distribution of $(\xi_{ij'}, j' \in [K])$, given that $\xi_{ij}$ contains a vertex with label $v'$. Assign label $v'$ to a uniformly chosen point in $\xi_{v''}^{j}$, and assign uniform labels without replacement from $V_j \setminus v'$ to the other points in $\xi_{v''}^{j}$ and uniform labels from $V_{j'}$ to the points in $\xi_{v''}^{j'}$ for $j' \in [K]\setminus j$.

If none of the newly assigned labels correspond to vertices in $\hat{V}^{e}(\ell)$ then ${E}^*(\ell+1,\ell')$ contains all edges in ${E}^*(\ell+1,\ell'-1)$ plus the edges with tail $v''$ and heads  corresponding to the labels of the points in $(\xi_{v''}^{j'}, j' \in [K])$, with the obvious edge lengths. In addition, all heads of those edges which were not in ${V}^{*,a}(\ell+1,\ell'-1)$ move to ${V}^{*,p}(\ell+1,\ell')$ (if they were not in that set already). Furthermore, $v'' \in {V}^{*,a}(\ell+1,\ell')$.

If any of the newly assigned labels correspond to vertices in $\hat{V}^{e}(\ell)$ then we flag the process and we return to the start of step $4)$. This last part of step 4b) is equivalent to conditioning on the event that there are no edges with tail $v''$ and an already explored vertex as head.

\item[4c)] Set $\hat G(\ell+1)=\hat G\left(\ell+1, \sum_{i=1}^Kx(v',i)\right)$.
\end{itemize}

\item[5)] Continue this process by increasing $\ell$ until there are no active vertices having a path of length less than $t^*$ towards $v$ in $\hat{G}(\ell)$. Say that $\ell^*$ is the smallest $\ell$ for which this is the case.
\end{itemize}

We note the following:
\begin{itemize}
\item The edge set $\hat{E}(\ell)$ is a subset of $E'$ and contains all edges in $\hat{G}(\ell)$. 
\item All edges in $E'$ with tails in $\hat{V}^{a}(\ell)$ are in $\hat{E}(\ell)$, but there might still be edges in $E'$ with heads in $\hat{V}^{a}(\ell)$  that are not in $\hat{E}(\ell))$. 
\item Every passive vertex is the head of an edge in $\hat{E}(\ell)$, but the passive vertices themselves are not explored, and none of the edges in $E'$ with tail in $\hat{V}^{p}(\ell)$ are in $\hat{E}(\ell)$. 
\item All edges in $E'$ with head or tail in $\hat{V}^{e}(\ell)$ are in $\hat{E}(\ell)$. 
\item The tails of edges in $\hat{E}(\ell)$ are in $\hat{V}^{a}(\ell) \cup \hat{V}^{e}(\ell)$ and their heads are in  $\hat{V}^{a}(\ell) \cup \hat{V}^{p}(\ell) \cup \hat{V}^{e}(\ell)$.
\end{itemize}

The construction of the exploration process yields $\hat{V}^e(\ell^*) = \mathcal{S}^{(n)}_v(t^*)$. Furthermore, if the process $\{\hat G(\ell)\}$ is not flagged until step $\ell^*$, then the construction of $\hat{V}^e(\ell^*)$ (and the distances from the vertices in this set to $v$) is equivalent to constructing a branching process $\{\mathcal{\tilde Z}^{(n)}(t);t \geq 0\}$ up to time $t^*$. In this branching process, particles of type $j$ give birth to a binomial distributed number of particles of type $i$, $i \in [K]$, where the parameters of the binomial distributed random variable are $n_i$ and $(n_j)^{-1}m_{ij}$. The number of children of the different types of particles are independent. Furthermore, the ages of the mother particles of type $j$ at birth of a child of type $i$ are independent and have density $\mathbb{E}[\xi_{ij}(da)]/m_{ij}$.


\subsection{The (backward) branching process}
\label{subsecbp}
We create a multi-type branching process $\{\mathcal{Z}(t);t \geq 0\}$ that can be coupled to $\{\mathcal{S}^{(n)}_v(t);t \geq 0\}$. The coupling is performed in Section~\ref{sec:coupling}. The branching process $\{\mathcal{Z}(t);t \geq 0\}$ is constructed in such a way that the distribution of the corresponding tree-like graph is the same as that of $\{\mathcal{S}^{(n)}_v(t);t \geq 0\}$ up to time $t^*$ with $t^*$ defined by~\eqref{eq:t*}. We leave out some of the details in the arguments. Those details can be filled in analogous to \cite{Ball95, Ball14, Barb13} for related models.

Without loss of generality we assume that $v \in V_1$. The multi-type (backward) branching process is as follows. The (single) ancestor is of type $1$. Particles of type $j$ give birth to particles of type $i$ according to a Poisson process with intensity $\frac{p_i}{p_j} \mathbb{E}[\xi_{ij}(da)]$, where $a$ is the age of the type $j$ particle, $i,j\in[K]$. All Poisson processes are independent.

The branching process $\{\mathcal{Z}(t);t \geq 0\}$ is analysed using existing theory from \cite{Iksa15}. First of all, the mean offspring measure of this backward branching process is defined through 
\begin{equation}
\label{meanmeas}
\mu^{(b)}_{ji}(dt)= \frac{p_i}{p_j} \mathbb{E}[\xi_{ij}(dt)].
\end{equation}
Note that 
\begin{equation*}
m^{(b)}_{ji}= \int_0^{\infty} \mu^{(b)}_{ji}(dt)= \int_0^{\infty} \frac{p_i}{p_j} \mathbb{E}[\xi_{ij}(dt)]=  \frac{p_i}{p_j}m_{ij}
\end{equation*} 
is the expected number of children of type $i$ of a particle of type $j$. Here $m_{ij}$ is given by~\eqref{eq:mij} (and the corresponding matrix is $M$). Let $M^{(b)}=\{m^{(b)}_{ji}\}_{j,i \in [K]}$. Straightforward matrix theory gives that $M$ and $M^{(b)}$ have the same dominant eigenvalue $R_0$, which by assumption is strictly larger than 1, i.e.\ the branching process is supercritical). 
Define 
\begin{equation}
\label{laplback}
\hat{m}^{(b)}_{ji}(x)= \int_{0}^{\infty} e^{-x t} \frac{p_i}{p_j} \mathbb{E}[\xi_{ij}(dt)]
\end{equation}
and let $\hat{M}^{(b)}(x)$ be the matrix with elements $\hat{m}^{(b)}_{ji}(x)$. Finally, let $\alpha$ be such that 
\begin{equation}
\label{maltdef}
\hat{M}^{(b)}_{ji}(\alpha)=1.
\end{equation}
Note that $K<\infty$, all elements of $M^{(b)}$ are finite and the dominant eigenvalue $R_0$ of $M^{(b)}$ is real and larger than 1. Therefore, $\alpha$ exists and is positive. 

We define the random vector $Z^i(t)$ as $Z^i(t) = (Z_1^i(t),Z_2^i(t),\cdots,Z_K^i(t))$, where $Z_j^i(t)$ is the number of particles of type $j \in [K]$ in $\mathcal{Z}(t)$ if the process starts with one newborn particle of type $i$. With some abuse of notation, let $\sigma(x)$ be the time of birth of particle $x$ in the branching process. Note that particles in the branching process $\mathcal{Z}(t)$ never die.  

We know that there exists an $\alpha >0$ such that 
\begin{equation}\label{bpeq1}
e^{-\alpha t} Z^{i}(t) \to W^{i} =(W^{i}_1, W^{i}_2, \cdots W^{i}_K) \qquad \mbox{a.s.\ as $t\to \infty$},
\end{equation}
where $W^{i}$ is a random vector that has, with probability 1, strictly positive elements on the set $\sum_{j=1}^K Z^{i}_j(t) \to \infty$ as $t \to \infty$~\cite{Iksa15}.

For particle $x \in \{\mathcal{Z}(t);t \geq 0\}$, define 
\begin{equation*}
\varphi_x(t) = \varphi_x(t;a) = \ind(t-\sigma(x)< a).
\end{equation*}
Let $\hat{Z}_j^{i}(t; a) = \sum_x  \varphi_x(t;a)$, where the sum is taken over all particles of type $j$ in $\{\mathcal{Z}(t);t \geq 0\}$, i.e.\ $\hat{Z}_j^{i}(t; a)$ is the number of particles of type $j$ that have age less than $a$ in the branching process at time $t$. Note that $\hat{Z}_j^{i}(t; a)$ is increasing in $a$.

If $Z_j^{i}(t) \to \infty$ for all $j \in [K]$ (i.e.\ if new particles keep on being born in the branching process), then
\begin{equation}
\label{bpeq2}
\frac{\hat{Z}_j^{i}(t; a)}{\sum_{j=1}^K Z_j^i(t)} \to c(a,j) \qquad \mbox{a.s.\ as $t \to \infty$,}
\end{equation}
where $c(a,j)$ is a constant independent of $i$ (Theorem 2.7 of \cite{Iksa15}). For our purposes we do not need to specify $c(a,j)$ further. 

Although the branching process is not dependent on $n$, we want to have some bound on the number of vertices born in the branching process as a function of $n$. This is used in the coupling in Section~\ref{sec:coupling}. We set 
\begin{equation}\label{eq:t*}
t^*= \frac{1-\kappa}{4}\log[n]/\alpha
\end{equation} 
(where $\kappa$ is as in Assumption \ref{boundass}). We obtain by (\ref{bpeq1}) that, for all $i,j \in [K]$,  
\begin{equation}
\label{martbound}
e^{-\alpha t^*}Z^{i}_j(t^*) = n^{-(1-\kappa)/4} Z^{i}_j(t^*) \to W^i_j \in (0,\infty) \qquad \mbox{ a.s.\ as $n \to \infty$}
\end{equation}
on the survival set of $\{\mathcal{Z}(t);t \geq 0\}$. In particular, this implies that
\begin{equation}
\label{martbound2}
 Z^{i}_j(t^*)  =o\left(n^{(1-\kappa)/3}\right) \qquad \mbox{w.h.p.}
\end{equation}

\subsection{The coupling}\label{sec:coupling}
Note that the only difference between the branching process $\{\mathcal{\tilde Z}^{(n)}(t)\}$ associated to the susceptibility process of Section~\ref{sec:susceptibility} and the multi-type branching process $\{\mathcal{Z}(t)\}$ of Section~\ref{subsecbp} is the distribution of the number of particles of type $i$  a particle of type $j$ gives birth to (binomially distributed with parameters $n_i$ and $m_{ij}/n_j$ and Poisson distributed with expectation $p_j/p_im_{ij}$, resp.)



We know from \cite[eq.\ (1.23)]{Barb92} (see also \cite{Barb13}) that the total variation distance between a binomial distributed random variable with parameters $n_i$ and $\frac{1}{n_i}\frac{n_i}{n_j} m_{ij}$ and a Poisson random variable with parameter $\frac{n_i}{n_j}  m_{ij}$ is $O(\frac{1}{n_i})$. Moreover, the total variation distance between a Poisson distributed random variable with parameter $\frac{n_i}{n_j}  m_{ij}$ and a Poisson distributed random variable with parameter $\frac{p_i}{p_j} m_{ij}$ is $O(\sqrt{|\frac{n_i}{n_j}-\frac{p_i}{p_j}|}) =  O(\frac{1}{\sqrt{n}})$ \cite[Theorem 1.C]{Barb92} (see also \cite{Barb13}). Here we have also used Assumption \ref{nipiass}. By the triangle inequality this implies that the total variation distance between a Poisson distributed random variable with expectation $\frac{p_i}{p_j} m_{ij}$ and a binomial distributed random variable with parameters $n_i$ and $\frac{1}{n_i}\frac{n_i}{n_j}  m_{ij}$ is $O(\frac{1}{\sqrt{n}})$. Hence, as long as the number of particles born in any of the two branching processes is $o(\sqrt{n})$, the two processes can be perfectly coupled w.h.p. In the remainder of this section we show that the coupling is w.h.p.\ perfect up to time $t^*$ with $t^*$ given by~\eqref{eq:t*}. 


By construction, if the first $|\hat{V}^a(\ell^*)| + |\hat{V}^e(\ell^*)|$ vertices that we ``try to include'' in  $\{\hat{V}^{a}(\ell) \cup \hat{V}^{e}(\ell);\ell \in \mathbb{N}_0\}$  are all different and none of the first $|\hat{V}^p(\ell^*)|$ that we ``try to include'' in  $\{\hat{V}^{p}(\ell);\ell \in \mathbb{N}_0\}$ are in  $\hat{V}^{a}(\ell^*) \cup \hat{V}^{e}(\ell^*)$, then the process is not flagged. 
 
The law of large numbers yields 
\begin{equation}
\label{setbound}
|\hat{V}^a(\ell^*)| \leq 2 \max_{i,j \in [K]}\frac{n_i}{n_j}m_{ij}  |\hat{V}^e(\ell^*)|, \qquad \mbox{w.h.p.}\end{equation}
Here we use that the expected number of edges in $E'$ with any given vertex in $V'$ as head is bounded from above by $\max_{i,j \in [K]}\frac{n_i}{n_j}m_{ij}$. Equation (\ref{setbound}) implies that 
\begin{equation}
\label{setbound1a}
|\hat{V}^e(\ell^*)| + |\hat{V}^a(\ell^*)|  = O(|\hat{V}^e(\ell^*)|), \qquad \mbox{w.h.p.}
\end{equation}

Using Assumption~\ref{boundass} we obtain
\begin{equation}
\label{setbound2}
|\hat{V}^p(\ell^*)| \leq n^{\kappa} (|\hat{V}^a(\ell^*)|+ |\hat{V}^e(\ell^*)|), \qquad \mbox{w.h.p.}
\end{equation}
Combining inequalities (\ref{setbound1a}) and (\ref{setbound2}) yields
\begin{equation}
\label{setbound3}
|\hat{V}^p(\ell^*)| = O_p(n^{\kappa}|\hat{V}^e(\ell^*)|), 
\end{equation}

Since $\{\mathcal{Z}(t);t \geq 0\}$ can be perfectly coupled with $\{\mathcal{S}^{(n)}_v(t);t \geq 0\}$ w.h.p.\ until $o(\sqrt{n})$ particles are born in $\{\mathcal{Z}(t);t \geq 0\}$, we obtain by using~\eqref{martbound2},  that $|\mathcal{Z}(t^*)| = o_p(n^{(1-\kappa)/3})$. Therefore, we also know that $|\mathcal{Z}(t^*)| = o(n^{(1/2})$ and $|\hat{V}^e(\ell^*)| = o_p(n^{(1-\kappa)/3})$.
Combined with (\ref{setbound3}) this gives
\begin{equation}
\label{setbound4}
|\hat{V}^e(\ell^*)| + |\hat{V}^a(\ell^*)| = o_p(n^{(1-\kappa)/3})
\end{equation}
and 
\begin{equation}
\label{setbound5}
|\hat{V}^p(\ell^*)| = o(n^{(1+2\kappa)/3}) = o_p(n).
\end{equation}
Using birthday-problem-like arguments \cite[p.24]{Grim01} the probability that the first $|\hat{V}^e(\ell^*)| + |\hat{V}^a(\ell^*)|$ activated vertices in $\{\hat{G}(\ell); \ell \in \mathbb{N}_0\}$ are not all different is bounded from above by 
\begin{equation*}
\left(\min_{j \in [K]}n_j\right)^{-1}(|\hat{V}^e(\ell^*)| + |\hat{V}^a(\ell^*)|)^2= o_p(n^{-(1+2\kappa)/3}),
\end{equation*}
while the probability that among the first $|\hat{V}^p(\ell^*)|$ vertices that we ``try to include'' in  $\{\hat{V}^{p}(\ell);\ell \in \mathbb{N}_0\}$ there are vertices in  $\hat{V}^{a}(\ell^*) \cup \hat{V}^{e}(\ell^*)$ is bounded from above by 
\begin{align*}
|\hat{V}^p(\ell^*)| \frac{|\hat{V}^e(\ell^*)| + |\hat{V}^a(\ell^*)|}{\min_{j \in [K]}n_j} &= o_p(n^{(1+2\kappa)/3}) \times o_p(n^{(1-\kappa)/3}) \times O(n^{-1}) \\
&=o_p(n^{(\kappa-1)/3})\\
&=o_p(1).
\end{align*}
We conclude that the probability that the exploration process  $\{\hat{G}(\ell); \ell \in \mathbb{N}_0\}$ is flagged up to and including step $\ell^*$ goes to 0 as $n \to \infty$, as we desired. We summarise the main result of this section in the following lemma.
\begin{lemma}
\label{couplelemma}
There exists a probability space on which we can define the branching process $\{\mathcal{Z}(t);t \geq 0\}$ and the susceptibility process $\{\mathcal{S}^{(n)}_v(t);t \geq 0\}$ for all $n$, such that
\begin{equation*}
\mathbb{P}\left(|\mathcal{S}^{(n)}_{v}(t;a,j)| = \hat{Z}_j^i(t; a)  \mbox{ for all $t \leq t^*$, $a \in (0,t^*)$ and $j \in [K]$}\right) \to 1,
\end{equation*}
where $t^*= t^*(n)$ is given by~\eqref{eq:t*}.
\end{lemma}

\section{Proofs}\label{sec:infector}

\subsection{Proof of Theorem~\ref{firstmain}}
We are interested in the expected fraction of vertices infected by vertices of type $i$ among the ultimately infected vertices of type $j$ during a major outbreak. By exchangeability this expected fraction is given by 
\begin{equation}
\label{infprob}
\mathbb{P}\left(\mbox{$v$ is infected by a type $i$ vertex}|\mbox{$v$ is ultimately infected},v \in V_j\right).
\end{equation}


Consider $\mathcal{S}^{(n)}_v(t^*)$, where we assign to each vertex $v' \in \mathcal{S}^{(n)}_v(t^*)$ the value $\sigma'_{v'} = d(v',v)$. In Section~\ref{sec:susceptibility} we constructed part of the graph $G'$ by looking backward in time. To analyse (\ref{infprob}), we construct another part of the graph $G'$, by looking forward in time and starting at $V_{\text{init}}$. By Assumption \ref{assvinit}, $\mathcal{S}^{(n)}_v(t^*)$ does not overlap with $V_{\text{init}}$ w.h.p. We condition on this event.

We construct the relevant part of $G'$ through a series of subgraphs $\{\tilde{G}(\ell); \ell \in \mathbb{N}_0\}$ as follows.

\begin{itemize}
\item[1)] Construct $\mathcal{S}^{(n)}_v(t^*)$ and the edges that connect these vertices in $G'$ as in Section~\ref{sec:susceptibility}. Let $\tilde{G}(0)$ be this graph together with (the isolated) vertices in $V_{\text{init}}$. 

\item[2)] The vertices in $V_{\text{init}}$ are active in $\tilde{G}(0)$. 

\item[3)] Assume that we know $\tilde{G}(\ell)$. Let $\tilde{\sigma}_v(\ell)$ be the distance from $V_{\text{init}}$ to $v$ in $\tilde{G}(\ell)$. If there is no path from $V_{\text{init}}$ to $v$ in $\tilde{G}(\ell)$ then we set $\tilde{\sigma}_v(\ell) = \infty$. If there are no active vertices in $\tilde{G}(\ell)$ with distance from $V_{\text{init}}$ less than $\tilde{\sigma}_v(\ell)-t^*$, then the shortest path from $V_{\text{init}}$ to $v$ is the same in $G'$ as in $\tilde{G}(\ell)$, and we set $\tilde{G}(k) = \tilde{G}(\ell)$, for all $k \geq \ell$.

\item[4)] If there are active vertices in $\tilde{G}(\ell)$ with distance from $V_{\text{init}}$ less than $\tilde{\sigma}_v(\ell)-t^*$, then we construct $\tilde{G}(\ell+1)$ as follows. Pick the active vertex with the lowest distance from $V_{\text{init}}$ in $\tilde{G}(\ell)$ (in case of a tie, which occurs if there are still several active vertices in $V_{\text{init}}$, make a uniform choice among those vertices). Say that this vertex is $u \in V_i$. Then assign to $u$ the point processes $(\xi_u^1, \cdots \xi_u^K)$, having the correct distribution (see Section \ref{secmod}) and label the points of $\xi_u^j$ ($j\in [K]$) with independent uniform vertices from $V_j$ without replacement.

First we check whether some of the labels chosen correspond to vertices in $\mathcal{S}^{(n)}_v(t^*)$. If there are such labels, we add ``preliminary'' edges with tail $u$ and heads equal to the respective vertices in $\mathcal{S}^{(n)}_v(t^*)$ (whether the preliminary edges become ``actual'' edges in the graph is determined in step 5 of the construction). The lengths of such edges correspond to the points in the point processes. 
Say that $u' \in \mathcal{S}^{(n)}_v(t^*) \cap V_j$ is one of the labels chosen and that $u$ is chosen at ``age'' $t'$. 

If $u' \not\in \mathcal{S}^{(n)}_v(t^*;t',j)$ for at least one of the vertices that $u$ connects to, then the distance from $u$ to $v$ becomes less than $t^*$ in the preliminary graph. We set $\tilde{G}(\ell+1) = \tilde{G}(\ell)$. Note that  we have identified all vertices from which there is a path of length at most $t^*$ to $v$ in the exploration of $\mathcal{S}^{(n)}_v(t^*)$. However, since we know that $u \not\in \mathcal{S}^{(n)}_v(t^*)$, we have to condition on the event that no edges with tail $u$ and heads in the set $\mathcal{S}^{(n)}_v(t^*-t')$ is shorter than $t'$.

\item[5)] Finally, in deciding whether the preliminary edges become edges in the graph $\tilde G(\ell)$, we consider the following. If the distance from $u$ to $v$ does not become less than $t^*$ in the preliminary graph, then $u$ becomes passive in $\tilde{G}(\ell+1)$, and the edge set of $\tilde{G}(\ell+1)$ consists of all edges in $\tilde{G}(\ell)$ plus the edges with tail $u$ and heads corresponding to the labels of the points in $(\xi_u^1, \cdots \xi_u^K)$. The edges that are added have the obvious lengths. The heads of the edges in $\tilde{G}(\ell+1)$ are part of the vertex set of $\tilde{G}(\ell+1)$. The active vertices in $\tilde{G}(\ell+1)$ are the active vertices in $\tilde{G}(\ell)$ apart from $u$ plus the newly added vertices in $\tilde{G}(\ell+1)$.

\end{itemize}

In the construction above, vertices of type $j$ are added in an interchangeable way, $j\in[K]$. If a preliminary edge of length $a$ has a head in $\mathcal{S}^{(n)}_v(t^*)\cap V_j$, then every vertex in $\mathcal{S}^{(n)}_v(t^*;a,j)$ has the same probability of being the head of this edge and becoming part of the constructed graph.

Next, for the coupling, we analyse $\mathcal{S}^{(n)}_v(t^*;a,j)$ even further. Denote the vertices in $\{u \in V_i; (u,v) \in E'\}$, i.e.\ the vertices in $V_i$, which are in the first generation of the susceptibility set of $v$, by $u_{i1},u_{i2}, \cdots, u_{i, X_{i}(v)}$. Here $$X_{i}(v) = |\{u \in V_i; (u,v) \in E'\}|$$ is the number of vertices of type $i$ in this first generation of the susceptibility set. Note that $X_{i}(v)$ is distributed as $X_{ij}$ where $X_{ij}$ is
binomially distributed with parameters $n_i$ and $m_{ij}/p_j$. So, $X_{ij}$ converges in distribution to a
Poisson distributed random variable with expectation $(p_i/p_j)m_{ij}$. For $r \in [X_{i}(v)]$, let $\tau_{ir}(v)= \eta(u_{ir},v)$ be the length of edge $(u_{ir},v) \in E'$. Here $\tau_{ir}(v)$ is distributed as $\tau_{ijr}$ with distribution function $\mathbb P(\tau_{ijr}\leq a)=\mathbb E[\xi_{ij}(a)]/m_{ij}$. For convenience, we let $X_i=X_i(v)$ and $\tau_{ir}=\tau_{ir}(v)$.

We consider the susceptibility processes of the vertices in $\{u_{ir}; i \in [K], r \in [X_{i}]\}$ up to distance $t^*$ from $v$ separately, i.e.\ we consider the susceptibility set of $u_{ir}$ up to time $t^*- \tau_{ir}$. Those susceptibility sets are (following the arguments in Section~\ref{subsecbp}) w.h.p.\  not overlapping and independent. Also note that
\begin{equation}
\label{susuni}
\mathcal{S}^{(n)}_{v}(t^*;a,j)= v \cup  \bigcup_{i =1}^K \bigcup_{r=1}^{X_{i}} \mathcal{S}^{(n)}_{u_{ir}}(t^*-\tau_{ir};a,j).
\end{equation}
By Lemma~\ref{couplelemma}, we know that 
\begin{equation}
\label{whpeq}
|\mathcal{S}^{(n)}_{u_{ir}}(t^*-\tau_{ir},a,j)| = \hat{Z}^i_j(t^*-\tau_{ir}; a) \qquad \mbox{w.h.p.}
\end{equation}
Furthermore, by equations (\ref{bpeq1}) and  (\ref{bpeq2}), we obtain that 
\begin{equation}
\label{infectconv}
e^{-\alpha t} \hat{Z}^{i}_j(t;  a) \to c(a,j) W^{i} \qquad \mbox{a.s.\ as $t \to \infty$}
\end{equation}
on $\{\hat{Z}^{i}(t) \to \infty\}$. Since the sets $\mathcal{S}^{(n)}_{u_{ir}}(t^*-\tau_{ir})$ are w.h.p.\ not overlapping (and do not contain $v$), it follows from (\ref{susuni}) that
\begin{equation*}
|\mathcal{S}^{(n)}_{v}(t^*;a,j)\setminus v|= \sum_{i =1}^K \sum_{r=1}^{X_{i}} |\mathcal{S}^{(n)}_{u_{ir}}(t^*-\tau_{ir};a,j)| \qquad \mbox{w.h.p.}
\end{equation*}
Then, by (\ref{whpeq}) and (\ref{infectconv}), for all $i \in [K]$ and $r \in [m_i]$,
\begin{equation*}
e^{-\alpha t^*} |\mathcal{S}^{(n)}_{u_{ir}}(t^*-\tau_{ir};a,j)| \to e^{-\alpha \tau_{ir}}  c(a,j) W^{i}(r)
\end{equation*}
in probability on  $\{\hat{Z}^{i}_j(t) \to \infty\}$ as $n \to \infty$. Here, the random variables $W^i(1)$, $W^i(2)$, $\ldots$ are independent copies of $W^i$, $i\in[K]$. Let $\hat{W}^i$ be distributed as $W^i\ind(\hat{Z}^{i}(t) \to \infty)$, and let $\hat{W}^i(r)$ be defined analogously. Note that $\{\hat{Z}^{i}(t) \to \infty\}$ implies that the branching process survives and new particles are born w.h.p.\ in the interval $[t^*/2,t^*]$. If $\{\hat{Z}^{i}(t) \not\to \infty\}$, then there is a last birth in the process. Because $t^*/2 \to \infty$ as $n \to \infty$, there is then no particle born in in the interval $[t^*/2,t^*]$ w.h.p. It follows from the coupling arguments above that, as $n \to \infty$,
\begin{equation*}
\mathbb{P}\left(\mathcal{S}^{(n)}_{v}(t^*) \neq \mathcal{S}^{(n)}_{v}(t^*/2)\right) \to \mathbb{P}\left(\hat{Z}^{i}(t) \to \infty\right).
\end{equation*}

We are interested in the fraction of vertices (possibly specified by type and age) at time $t^*$ in $\mathcal{S}^{(n)}_v(t^*)$ that are connected to $v$ through a path of vertices in $\mathcal{S}^{(n)}_v(t^*)$ that include $u_{ir}$, $i \in [K]$ and $r \in [X_{ij}]$. That is, we want to analyse 
\begin{equation*}
\frac{|\mathcal{S}^{(n)}_{u_{ir}}(t^*-\tau_{1r};a,j)|}{\sum_{j =1}^K \sum_{r=1}^{X_{j}} |\mathcal{S}^{(n)}_{u_{jr}}(t^*-\tau_{jr};a,j)|},
\end{equation*}
on the set $\sum_{j =1}^K \sum_{r=1}^{X_{j}} |\mathcal{S}^{(n)}_{u_{jr}}(t^*-\tau_{jr};a,j)| \neq 0$. Lemma \ref{couplelemma} allows us to couple the epidemic process with a branching processes $\{\mathcal{Z}(t);t \geq 0\}$ such that w.h.p.
\begin{equation}\label{eq:fraction}
\frac{|\mathcal{S}^{(n)}_{u_{ir}}(t^*-\tau_{1r};a,j)|}{\sum_{j =1}^K \sum_{r=1}^{X_{j}} |\mathcal{S}^{(n)}_{u_{jr}}(t^*-\tau_{jr};a,j)|} = 
\frac{Z_j^{i,r}(t^*-\tau_{1r};a)}{\sum_{k =1}^K \sum_{r=1}^{X_{kj}} Z_j^{k,r}(t^*-\tau_{kr};a)},
\end{equation}
where $Z_j^{i,r}(t;a)$ are independent copies of $Z_j^{i}(t;a)$, $r \in \mathbb{N}$. Multiplying numerator and denominator of~\eqref{eq:fraction} by $e^{-\alpha t^*}$ and using (\ref{infectconv}) we obtain that (\ref{eq:fraction}) is equal to
\begin{equation*}
\frac{e^{-\alpha \tau_{ir}} (e^{-\alpha (t^*-\tau_{ir})}Z_j^{i,r}(t^*-\tau_{1r};a))}{\sum_{k =1}^K \sum_{r=1}^{X_{kj}} e^{-\alpha \tau_{kr}} (e^{-\alpha (t^*-\tau_{kr})} Z_j^{k,r}(t^*-\tau_{kr};a)},
\end{equation*}
which converges a.s.\ to
\begin{equation}\label{eq:fraction2}
\frac{e^{-\alpha \tau_{ir}} c(a,j) \hat{W}^i(r)}{\sum_{k =1}^K \sum_{r=1}^{X_{kj}} e^{-\alpha \tau_{kr}} c(a,j) \hat{W}^k(r)} = \frac{e^{-\alpha \tau_{ir}}  \hat{W}^i(r)}{\sum_{k =1}^K \sum_{r=1}^{X_{kj}} e^{-\alpha \tau_{kr}}  \hat{W}^k(r)}
\end{equation}
as $n \to \infty$ (i.e.\ as $t^* \to \infty$). Note that the right hand side of~\eqref{eq:fraction2} does not depend on $a$ or $j$. Since vertices of type $j$, $j \in [K]$, are added in such a way that every vertex in $\mathcal{S}^{(n)}_v(t^*;a,j)$ has the same probability of being the head of this edge, every vertex in $\mathcal{S}^{(n)}_v(t^*;a,j)$ has the same probability of becoming part of the constructed graph. So, 
$\mathbb{P}\left(\mbox{$v$ is infected by a type $i$ vertex}|\mbox{$v$ is ultimately infected}\right)$ is equal to
$$
\frac{1}{\mathbb{P}(\mathcal{S}^{(n)}_{v}(t^*) \neq \mathcal{S}^{(n)}_{v}(\frac{t^*}{2}))} \mathbb{E}\left[\frac{\sum_{r=1}^{X_i} |\mathcal{S}^{(n)}_{u_{ir}}(t^*-\tau_{ir};a,j)|}{\sum_{k =1}^K \sum_{r=1}^{X_k} |\mathcal{S}^{(n)}_{u_{kr}}(t^*-\tau_{kr};a,j)|} \ind\left(\mathcal{S}^{(n)}_{v}(t^*) \neq \mathcal{S}^{(n)}_{v}\left(\frac{t^*}{2}\right)\right)
\right]
$$
and converges to $\rho_{ij}$, with $\rho_{ij}$ given by~\eqref{mainequ}. This completes the proof for Theorem \ref{firstmain}.



\subsection{Towards to proof of Theorem~\ref{secondmain}: Bounds for the $\boldsymbol{\rho_{ij}}$}\label{sec:bounds}
Theorem~\ref{firstmain} provides us with an expression for the asymptotic fractions $\rho_{ij}$ of infected individuals of type $j$ that were infected by individuals of type $i$. However, often there is no explicit description of the distribution of ${W}^k(r)$. In Theorem~\ref{secondmain} we consider bounds for the $\rho_{ij}$ for a special class of models, as specified in Section \ref{secmod}. In order to obtain those bounds, in this subsection and Section \ref{sec:proofsecondmain}, we discuss, for the general setting, how to obtain the bounds for $\rho_{ij}$ using the epidemic random graph $G'$. In most of the subsequent analysis we analyse the graph $G'$ without taking the lengths of edges in into account. 


\begin{figure}
\centering
\includegraphics[scale=0.7]{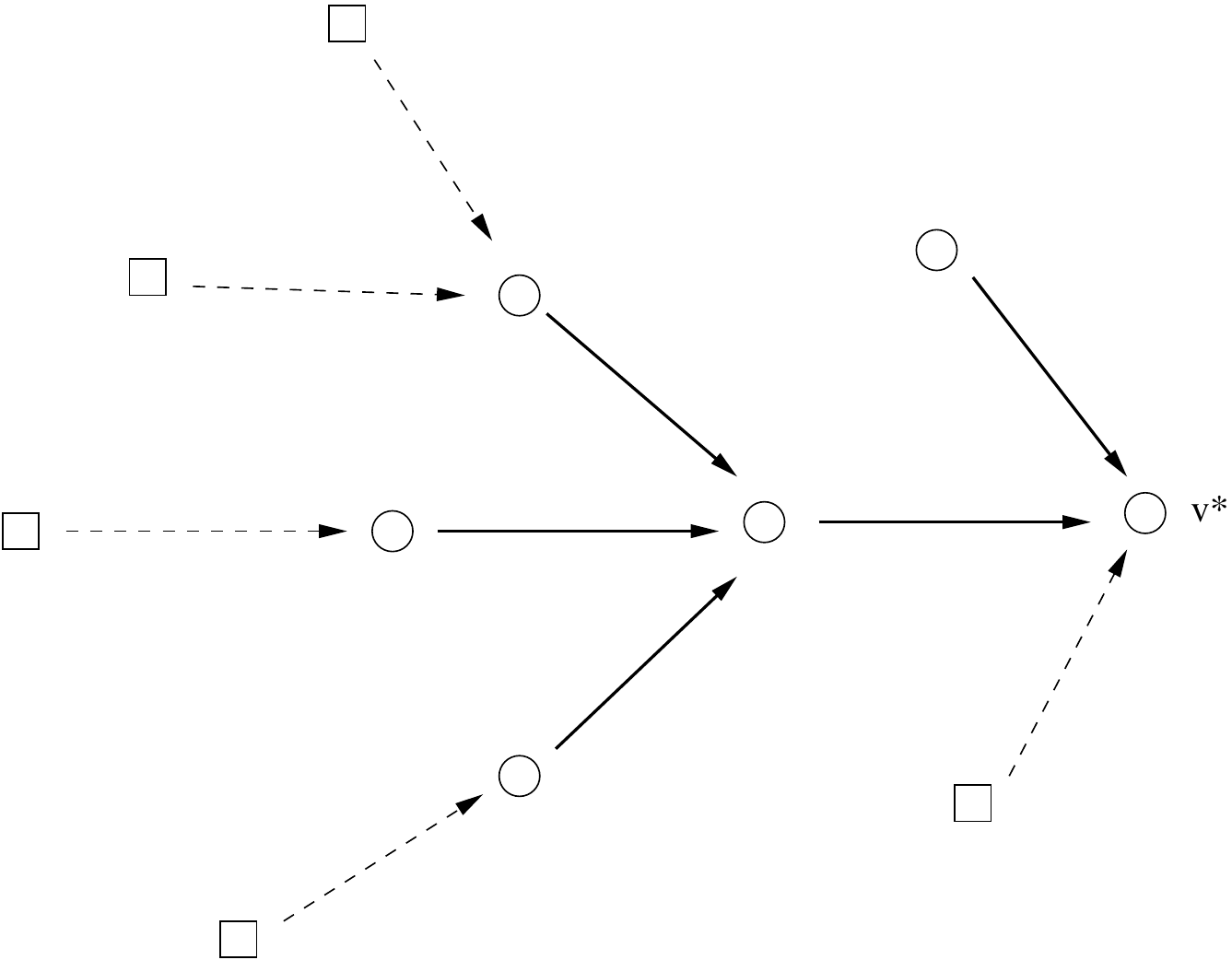}
\caption{Illustration of (a part of) the susceptibility set of $v_*$ in $G'=(V,E')$. Vertices of type 1 are represented by circles and vertices of type 2 are represented by boxes. The part of the susceptibility set illustrated in the figure is the set connected to $v_*$ through paths of edges with heads of type 1.}
\label{fig:backward}
\end{figure}

Note that Assumption \ref{finass} guarantees that, with probability 1, all paths in $E'$ have different lengths. In order to obtain the maximum and minimum of the probability $\rho_{ij}$ for fixed mean offspring matrix $M$, we first investigate the susceptibility set of $v_*$ ($v_* \in V_j$) restricted to the graph $G'_{ij}=(V,E'_{ij})$. We denote this susceptibility set by $\mathcal{S}_{v_*,ij}$. Here $E'_{ij} \subset E'$ is the subset of $E'$ that consists of the edges that either have tail vertex in $V_i$ or head vertex not in $V_j$, i.e.\ $E \setminus E'_{ij}$ is the set of edges with tails in $V \setminus V_i$ and heads in $V_j$.
In Figure \ref{fig:backward}, the set $\mathcal{S}_{v_*,11}$ consist of all  vertices of type 1 (the circles).

Let $\mathcal{S}^{j}_{v_*,ij}= V_j \cap \mathcal{S}_{v_*,ij}$. If $\mathcal{S}_{v_*,ij} \cap V_{\text{init}} = \emptyset$, i.e.\ if there is no path from $V_{\text{init}}$ to $v_*$ in $G'_{ij}$, then, by the definition of the epidemic process, every  vertex of type $j$ in $\mathcal{S}_{v_*,ij}$ has the same probability to be the first one to be infected in the epidemic, i.e.\ $u_*=\argmin_ {u \in \mathcal{S}^j_{v_*,ij}}d(V_{\text{init}},u)$ is uniform in $\mathcal{S}^j_{v_*,ij}$.

Condition on $\mathcal{S}_{v_*} \cap V_{\text{init}} \neq \emptyset$. If $u_* =v_*$ and $\mathcal{S}_{v_*,ij}\cap V_{\text{init}} = \emptyset$, then $v_*$ is infected by a vertex that is not of type $i$. On the other hand, if $u_* \neq v_*$ or $\mathcal{S}_{v_*,ij} \cap V_{\text{init}} \neq \emptyset$, then $v_*$ might be infected by a type $i$ vertex. Hence
\begin{align*}
1-\rho_{ij} &= \mathbb{P}(\mbox{$v_*$ is not infected by a type $i$ vertex}|\mathcal{S}_{v_*} \cap V_{\text{init}} \neq \emptyset)\\
& \geq \mathbb{P}(u_* = v_*|\mathcal{S}_{v_*} \cap V_{\text{init}} \neq \emptyset).
\end{align*}

Now assume that $\mathbb{P}(\eta_{i',j'}> n^{-1}) =1$, for $i' \in [K]\setminus i$ and $j'=j$, with $j$ the type of $v_*$, and $\mathbb{P}(\eta_{i',j}< n^{-2}) =1$ otherwise. This implies that the lengths of any path with only edges in $E'_{ij}$ is less than any edge in $E'\setminus E'_{ij}$. The assumptions guarantee that, for $u_* \neq v_*$ or $\mathcal{S}_{v_*,ij} \cap V_{\text{init}} \neq \emptyset$, conditioned on $\mathcal{S}_{v_*} \cap V_{\text{init}} \neq \emptyset$, $v_*$ is infected by a vertex of type $i$ (the tail of the edge with head $v_*$ in the shortest path in $E_{ij}'$ from $u_*$ or $V_{\text{init}}$ to $v_*$). Therefore, for this model,
\begin{equation*}
\mathbb{P}(\mbox{$v_*$ is not infected by a type $i$ vertex}|\mathcal{S}_{v_*} \cap V_{\text{init}} \neq \emptyset) = \mathbb{P}(u_* = v_*|\mathcal{S}_{v_*} \cap V_{\text{init}} \neq \emptyset).
\end{equation*}
Hence, for a given distribution of $E'$, models with $\mathbb{P}(\eta_{i',j'}> n^{-1}) =1$, for $i' \in [K]\setminus i$ and $j'=j$  and $\mathbb{P}(\eta_{i',j}< n^{-2}) =1$ otherwise, are among the models for which the fraction of the ultimately infected vertices of type $j$, infected by a type $i$ vertex is maximised. 

Using similar arguments we obtain that $\rho_{ij}$, the fraction of ultimately infected  vertices of type $j$ that are infected by vertices of type $i$, is minimal if the edge lengths of vertices with tail in $i$ and head in $j$ are much longer than the other edges. This will be used in Section \ref{sec:proofsecondmain}.

\subsection{Proof of Theorem~\ref{secondmain}}\label{sec:proofsecondmain}
We consider the so-called symptom-response SEIR epidemic model that is introduced in detail in \cite{Leun18}, see also Remark~\ref{rk:SEIR}.
That is, consider the model introduced in Section \ref{secmod} with $K=2$ and consider general distributions $(\xi_{i1}, \xi_{i2})$. 
For this model, using the arguments from Section \ref{sec:bounds}, we can compute the maximal probabilities $\rho_{11}^+$, $\rho_{22}^+$ (and minimal probabilities $\rho_{21}^-$ and $\rho_{11}^-$) explicitly. 
In this subsection we consider the model in which edges from $V_1$ to $V_1$ are infinitesimally short and all other edges in $E'$ are relatively long. 
From Section \ref{sec:bounds} we know that this is the model for which the fraction of vertices of type $1$ infected by type $1$ vertices is maximised. 
For reasons of convenience we assume that $v_*$ is of type 1 with $\mathcal S_{v_*}\cap V_{\text{init}}\neq\emptyset$. Note that Theorem~\ref{secondmain} considers the special case that $(\xi_{i1}, \xi_{i2})$ is obtained from the independent labelling of a one-dimensional point process $\xi_i$. In this special case we have $\rho_{11}=\rho_1=\rho_{12}$ and $\rho_{21}=\rho_2=\rho_{22}$. We treat Theorem~\ref{secondmain} in Remark~\ref{rk:prf_thm2} at the end of this section.
First, we compute the upper bound $\rho_{11}^+$ (and lower bound $\rho_{21}^-$) for the general setting. 
We note that it is harder (if not impossible) to obtain an explicit expression for $\rho_{12}^+$ or $\rho_{21}^+$. As will become clear in the computation below, the difficulty with $\rho_{12}^+$ or $\rho_{21}^+$ is that one would need to consider paths in $G_{12}'$ and $G_{21}'$ that contain vertices of both type 1 and type 2. In contrast, for computing $\rho_{11}^+$,  paths in $G_{11}'$ that end in $v^*$ contain only vertices of  type 1.


If we ignore the lengths of edges in $G'$ in the general model introduced in Section \ref{secmod}, then
the approximating (backward) branching process describing the generation-based growth of $\mathcal{S}_{v_*}(t)$ is defined through the following offspring distributions. The number of children of type $i$ of a particle of type $j$ is Poisson distributed with expectation $m_{ji}^{(b)} = \frac{p_i}{p_j}m_{ij}$, $i,j\in[2]$. For different $i$ and $j$ the distributions are independent of each other.

Let $Y= |\mathcal{S}_{v_*,11}|$. It is easily seen that $Y$ is approximated by the size of a branching process with Poisson offspring distribution that has expectation $m_{11}^{(b)}$. Then $Y$ is Borel distributed with parameter $m_{11}^{(b)}$, i.e.\ for $\ell \in \{1,2,\cdots\}$, 
\begin{equation}
\label{Boreldist}
\mathbb{P}(Y=\ell) = \frac{(m_{11}^{(b)} \ell)^{\ell-1} e^{-m_{11}^{(b)} \ell}}{\ell!}
\end{equation}
(see~\cite{Aldo98}). If $m_{11}^{(b)}>1$, then $\mathbb{P}(Y=\infty) > 0$. Standard results on Borel distributions \cite{Aldo98} give that, for $m_{11}^{(b)} \leq 1$,
\begin{equation}
\label{Borelinv}
 \mathbb{E}[1/Y] = 1-m_{11}^{(b)}/2.
\end{equation}
Define
\begin{equation*}
\rho_{11} = \mathbb{P}(\mbox{$v_*$ is infected by by a type 1 vertex}|\mathcal{S}_{v_*} \cap V_{\text{init}} \neq \emptyset),
\end{equation*}
for $v_* \in V_1$. From the arguments in Section~\ref{sec:bounds} we know that
\begin{equation}
\label{rho1eq}
\rho_{21}^- = \mathbb{E}[1/Y|\mathcal{S}_{v_*} \cap V_{\text{init}} \neq \emptyset.]
\end{equation}
Consistency then yields
\begin{equation}\label{eq:rho11}
\rho_{11}^+=1-\rho_{21}^-.
\end{equation}

We use the following extinction probabilities in the backward branching process:
\begin{itemize}
\item[$q_1$:] the probability that the backward branching process starting with a single type 1 particle goes extinct,
\item[$q_2$:] the probability that the backward branching process starting with a single type 2 particle goes extinct,
\item[$\tilde{q}_1$:] the probability that the backward branching process restricted to type 1 particles dies out, i.e.\ $\tilde{q}_1 = \mathbb{P}(Y<\infty)$.
\end{itemize}

From theory on multi-type supercritical branching processes \cite[Chap.\ 4]{Jage75} we know that $(q_1, q_2)$ is the unique solution in $(0,1)^2$ of
\begin{eqnarray}
x & = & \sum_{k=0}^{\infty} \frac{(m_{11}^{(b)})^k}{k!}e^{-m_{11}^{(b)}} x^k
\sum_{\ell=0}^{\infty} \frac{(m_{12}^{(b)})^{\ell}}{\ell!}e^{-m_{12}^{(b)}} y^{\ell} 
= e^{-[m_{11}^{(b)}(1-x)+ m_{12}^{(b)}(1-y)]} \label{q1eq}\\
y & = & \sum_{k=0}^{\infty} \frac{(m_{21}^{(b)})^k}{k!}e^{-m_{21}^{(b)}} x^k
\sum_{\ell=0}^{\infty} \frac{(m_{22}^{(b)})^{\ell}}{\ell!}e^{-m_{22}^{(b)}} y^{\ell} 
= e^{-[m_{21}^{(b)}(1-x)+ m_{22}^{(b)}(1-y)]}. \label{q2eq}
\end{eqnarray}
Furthermore, $\tilde{q}_1$ is the smallest positive solution of~\eqref{qsequa} which might be 1 or strictly smaller than 1, depending on whether or not $p_1m_{11} = m_{11}^{(b)} \leq 1$.

We distinguish between $m_{11}^{(b)} \leq 1$ and $m_{11}^{(b)}>1$. This distinction is not necessary, but we think the argument becomes clearer, by treating the case $\tilde{q}_1= 1$ separately.

\subsubsection{The case $\boldsymbol{m_{11}^{(b)}  \leq 1}$}
Assume that $m_{11}^{(b)}  \leq 1$. Let $\mathcal{A}$ be the event that the backward branching process with ancestor of type 1, involving both type 1 and type 2 individuals, survives. We explore the backward branching process of a particle of type 1 on $\mathcal{A}$ as follows.

\begin{itemize}
\item Explore the backward process of particles of type 1. If we ignore the conditioning on $\mathcal{A}$, the process can be described by a subcritical branching process with Poisson offspring distribution with expectation $m_{11}^{(b)}$. The random variable  $Y$ is the total size of this branching process, including the initial individual. We know that $Y$ is Borel($m_{11}^{(b)}$) distributed.

\item Condition on event $\mathcal{A}$: the probability that a particle of type 1 has infinitely many descendants is $\mathbb{P}(\mathcal{A})= 1-q_1$, where $q_1$ is defined through~\eqref{q1eq} and \eqref{q2eq}. We also use the probability that a particle has infinite offspring, conditioned on having no children of type 1. This probability is $1-q$, where
\begin{equation}\label{eq:qq1}
q= \sum_{\ell =0}^{\infty} \frac{(m_{12}^{(b)})^{\ell}}{\ell !} e^{-m_{12}^{(b)}} (q_2)^{\ell} = e^{-m_{12}^{(b)}(1-q_2)}=q_1 e^{m_{11}^{(b)} (1-q_1)}.
\end{equation}
Using Bayes' rule we obtain
\begin{align}
\mathbb{P}(Y=\ell|\mathcal{A}) & = \frac{\mathbb{P}(\mathcal{A}|Y=\ell) \mathbb{P}(Y= \ell)}{1-q_1}\nonumber\\
& = \frac{1-q^{\ell}}{1-q_1}  \frac{(m_{11}^{(b} \ell)^{\ell-1} e^{-m_{11}^{(b)} \ell}}{\ell!}\nonumber\\
&= \frac{1}{1-q_1}  \frac{(m_{11}^{(b)} \ell)^{\ell-1} e^{-m_{11}^{(b)} \ell}}{\ell!} - \frac{q}{1-q_1}  \frac{(m_{11}^{(b)} \ell q)^{\ell-1} e^{-m_{11}^{(b)} \ell}}{\ell!}\nonumber\\
&=\frac{1}{1-q_1} \left(  \frac{(m_{11}^{(b)} \ell)^{\ell-1} e^{-m_{11}^{(b)} \ell}}{\ell!} -  \frac{q_1(m_{11}^{(b)} q_1\ell)^{\ell-1} e^{-m_{11}^{(b)} q_1 \ell}}{\ell!} \right).\label{Borelcond}
\end{align}
where we used~\eqref{eq:qq1} in the last equality.

\item Using~\eqref{Borelinv}, \eqref{rho1eq}  and (\ref{Borelcond}) we find that 
\begin{align}
\label{hitprob}
\rho_{21}^-&= \mathbb{E}[Y^{-1}|\mathcal{A}]\nonumber\\
&= \frac{1}{1-q_1} \left(1- \frac{m_{11}^{(b)}}{2} -q_1\left(1-\frac{m_{11}^{(b)} q_1}{2} \right) \right)\nonumber\\
&= \frac{1}{1-q_1} \left(1- q_1 - (1-(q_1)^2) \frac{m_{11}^{(b)}}{2}  \right)\nonumber\\
&= 1 - (1+q_1) \frac{m_{11}^{(b)}}{2}.
\end{align}

\item Using~\eqref{eq:rho11} we find the desired expression $\rho_{11}^+$ for $m_{11}^{(b)}$:
\begin{equation}\label{eq:upper}
\rho_{11}^+=1-(1+q_1) \frac{m_{11}^{(b)}}{2}.
\end{equation}

\end{itemize}

\subsubsection{The case $\boldsymbol{m_{11}^{(b)}  > 1}$}
Assume that $m_{11}^{(b)}  > 1$. Then it is possible that the backward process restricted to particles of type 1 is already large, i.e.\ the approximating branching process with Poisson($m_{11}^{(b)}$) offspring distribution already survives (call this event $\mathcal{A}_1$). In that case the probability that vertex $v$ is infected by a vertex of type 1 approaches 1 as the population size tends to infinity. 

The other possibility is that the backward process restricted to particles of type 1 stays small. Call this event $\mathcal{A}_1^C$, the complement of $\mathcal{A}_1$. The probability of this event is $\mathbb{P}(\mathcal{A}_1^C)=\tilde{q}_1$, where $\tilde{q}_1$ is the unique solution in $(0,1)$ of equation (\ref{qsequa}) (note that $m_{11}^{(b)}>1$ ensures that $\tilde q_1$ exists). From the theory of branching processes we know that, conditioned on $\mathcal{A}_1^C$, the approximating branching process is still a branching process with Poisson distributed offspring distribution, but now with offspring expectation $m_{11}^{(b)} \tilde{q}_1$. We still condition on the event that the approximating branching process (including both types of particles) survives.

Conditioned on $\mathcal{A}_1^C$, the total size of the backward process restricted to particles of type 1 is Borel distributed with parameter $m_{11}^{(b)} \tilde{q}_1$. We use Bayes' rule, 
\begin{equation*}
\mathbb{P}(Y=\ell|\mathcal{A},\mathcal{A}_1^C) = \frac{\mathbb{P}(\mathcal{A}|Y=\ell,\mathcal{A}_1^C) \mathbb{P}(Y= \ell|\mathcal{A}_1^C)}{\mathbb{P}(\mathcal{A}|\mathcal{A}_1^C)}.
\end{equation*}
Note that 
\begin{equation*} 
\mathbb{P}(\mathcal{A}|\mathcal{A}_1^C) = \frac{ \mathbb{P}(\mathcal{A},\mathcal{A}_1^C)}{ \mathbb{P}(\mathcal{A}_1^C)}= \frac{\mathbb{P}(\mathcal{A}_1^C)-\mathbb{P}(\mathcal{A}^C,\mathcal{A}_1^C)}{ \mathbb{P}(\mathcal{A}_1^C)}= \frac{\mathbb{P}(\mathcal{A}_1^C)-\mathbb{P}(\mathcal{A}^C)}{ \mathbb{P}(\mathcal{A}_1^C)}= \frac{\tilde{q}_1-q_1}{\tilde{q}_1}.
\end{equation*}
Then, with the same arguments as those leading to (\ref{hitprob}) but with $m_{11}^{(b)}$ replaced by $m_{11}^{(b)} \tilde{q}_1$, yield
\begin{equation*}
\mathbb{P}(\mathcal{A}|Y=\ell,\mathcal{A}_S^C) =1- q^{\ell} = 1-(q_1)^{\ell} e^{m_{11}^{(b)} (1-q_1)\ell}.
\end{equation*}
and 
\begin{equation*}
\mathbb{P}(Y= \ell|\mathcal{A}_1^C)= 
\frac{(m_{11}^{(b)} \tilde{q}_1 \ell)^{\ell-1} e^{-m_{11}^{(b)} \tilde{q}_1 \ell}}{\ell!}.
\end{equation*}
Combining these identities with~\eqref{qsequa} yields
\begin{align*}
& \mathbb{P}(Y=\ell|\mathcal{A},\mathcal{A}_1^C)\\
 &= \frac{\tilde{q}_1}{\tilde{q}_1-q_1}  \left(1-(q_1)^{\ell} e^{m_{11}^{(b)} (1-q_1)\ell}\right)\left(\frac{(m_{11}^{(b)} \tilde{q}_1 \ell)^{\ell-1} e^{-m_{11}^{(b)} \tilde{q}_1 \ell}}{\ell!} \right)\\
 &=  \frac{\tilde{q}_1}{\tilde{q}_1-q_1}  \left(\frac{(m_{11}^{(b)} \tilde{q}_1 \ell)^{\ell-1} e^{-m_{11}^{(b)} \tilde{q}_1 \ell}}{\ell!} \right) - \frac{q_1 (\tilde{q}_1)^{\ell}}{\tilde{q}_1-q_1} \left(\frac{(m_{11}^{(b)} q_1 \ell)^{\ell-1} e^{-m_{11}^{(b)} (\tilde{q}_1+q_1-1) \ell}}{\ell!} \right)\\
& =  \frac{\tilde{q}_1}{\tilde{q}_1-q_1} \left(\frac{(m_{11}^{(b)} \tilde{q}_1 \ell)^{\ell-1} e^{-(m_{11}^{(b)} \tilde{q}_1 \ell}}{\ell!} \right) - \frac{\tilde{q}_1}{\tilde{q}_1-q_1}\left(\frac{(m_{11}^{(b)} q_1 \ell)^{\ell-1} e^{-m_{11}^{(b)} q_1 \ell}}{\ell!} \right).
\end{align*}
By standard results on Galton Watson branching processes we know that $m_{11}^{(b)} \tilde{q}_1 \leq 1$. Moreover, since $q_1 \leq \tilde{q}_1$, also $m_{11}^{(b)} q_1 \leq 1$ holds. Then, using (\ref{Borelinv}), we obtain
\begin{align*}
\mathbb{E}[Y^{-1}|\mathcal{A},\mathcal{A}_1^C] &= \frac{\tilde{q}_1}{\tilde{q}_1-q_1}  \left(1-\frac{m_{11}^{(b)} \tilde{q}_1}{2}\right) - 
 \frac{q_1}{\tilde{q}_1-q_1} \left(1-\frac{m_{11}^{(b)} q_1}{2}\right)\\
& = 1 -  \frac{m_{11}^{(b)} (\tilde{q}_1 + q_1)}{2}.
\end{align*}
This leads to 
\begin{align}
\rho_{21}^- &= \mathbb{E}[Y^{-1}|\mathcal{A}]\nonumber\\
&= \mathbb{E}[Y^{-1}|\mathcal{A},\mathcal{A}_1^C]\mathbb{P}[\mathcal{A}_1^C|\mathcal{A}] + \mathbb{E}[Y^{-1}|\mathcal{A},\mathcal{A}_1] \mathbb{P}[\mathcal{A}_1|\mathcal{A}]\nonumber\\
&= \mathbb{E}[Y^{-1}|\mathcal{A},\mathcal{A}_1^C]\frac{\mathbb{P}[\mathcal{A}_1^C,\mathcal{A}]}{\mathbb{P}[\mathcal{A}]} + 0 \nonumber\\
&= \mathbb{E}[Y^{-1}|\mathcal{A},\mathcal{A}_1^C]\frac{\mathbb{P}[\mathcal{A}_1^C]-\mathbb{P}[\mathcal{A}_1^C,\mathcal{A}^C]}{\mathbb{P}[\mathcal{A}]}\nonumber\\
&= \left( 1 -  \frac{m_{11}^{(b)} (\tilde{q}_1 + q_1)}{2} \right) \frac{\tilde{q}_1-q_1}{1-q_1}.\label{finaleq}
\end{align}
Note that this expression~\eqref{finaleq} is consistent with the result~\eqref{hitprob} for $m_{11}^{(b)}\leq 1$, where $\tilde{q}_1=1$. Finally, using~\eqref{eq:rho11}, we find the expression for $\rho_{11}^+$:
\begin{equation}\label{eq:rho11_final}
\rho_{11}^+=1-\left( 1 - \frac{m_{11}^{(b)} (\tilde{q}_1 + q_1)}{2} \right) \frac{\tilde{q}_1-q_1}{1-q_1}.
\end{equation}


\begin{remark}[Symptom-response SEIR epidemic model: proof of Theorem~\ref{secondmain}]\label{rk:prf_thm2}
In \cite{Leun18} we interpret vertices of type 1 as individuals that show symptoms when infectious, while vertices of type 2 are asymptomatic throughout their infectious period. We assume that whether infected individuals become symptomatic or not, does not depend on who infected them or when they were infected. So, we may assign i.i.d.\ types to the vertices before the epidemic. Note that the type of a vertex that does not become infected during the epidemic has no epidemiological relevance. Type 1 and type 2 vertices are equally susceptible. This implies that $\{(\xi_v^1(t),\xi_v^2(t)); t \geq 0\}$ can be obtained by considering a one-dimensional point process $\{\xi_v(t); t \geq 0\}$, for all $v\in V$.  One can assign the types (type 1 with probability $p_1$ or type 2 with probability $p_2=1-p_1$) independently to the points of this process. Then $\rho_{11}=\rho_1=\rho_{12}$ and $\rho_{22}=\rho_2=\rho_{21}$. Furthermore $m_{ji}^{(b)}=p_i \tilde m_j$, with $\tilde m_j=\mathbb{E}[\xi_j(\infty)]$ the expected number of secondary cases generated by a newly infected type $j$ individual in an otherwise susceptible population, $i,j=1,2$. Note that~\eqref{eq:rho11_final} reduces to the fraction $\rho_1^+$ of the ultimately infected vertices, that are infected by symptomatic vertices: 
\begin{equation*}
\rho_1^+ = 1-\left( 1 - \frac{p_1\tilde m_1(\tilde{q}_1 + q)}{2} \right) \frac{\tilde{q}_1-q}{1-q},
\end{equation*}
where $\tilde{q}_1$ and $q$ are solutions of~\eqref{qsequa} and~\eqref{qequa}. This proves Theorem~\ref{secondmain}.
\end{remark}

\section{Discussion}\label{sec:discussion}

In this manuscript we couple a (fairly) general multi-type stochastic epidemic process with a weighted random graph. We use this random graph to obtain a characterisation of the (large population limit) fraction of individuals in the population that had an infector of (say) type 1 given a large outbreak. 
The results of this paper are applied (and in more detail in \cite{Leun18}) to a model where the types of individuals represent whether the infected individual will show symptoms at some moment if infected, or he or she stays asymptomatic  if  infected. 

From a public health perspective, a relevant variant of the model would be to consider a population with only one type of individuals, where infectious individuals may start off asymptomatic after which they become symptomatic (see \cite{Fras04}, where the asymptomatic phase is also referred to as ``presymptomatic''). The related question for this model would be: ``what fraction of the infected population had a symptomatic infector, given a large outbreak?''. In order to analyse this model we need to assign types to the edges in the epidemic graph instead of to the vertices. We can use the techniques of this paper to characterise the answer to the above question. We can still define a susceptibility process from a vertex $v$. The susceptibility process can be approximated by a multi-type branching process. The type of a particle in the branching process should then corresponds to the type of the edge through which the particle is added to the susceptibility process. That is, the type of the particle (say $u$) in the branching process depends on whether the first edge in the shortest path from $u$ to $v$ in the epidemic graph  represents a contact that was made while $u$ was symptomatic or asymptomatic. As in Theorem \ref{firstmain} the answer will be implicit through its dependence on the distribution of the martingale limits $W$.

\section*{Acknowledgements}
The authors are supported by Vetenskapsr{\aa}det (Swedish Research Council), grants 2015-05015
(TB and KYL) and  2016-04566 (PT).

\bibliographystyle{abbrv}
\bibliography{publications.bib}

\end{document}